\documentclass[12pt]{article}
\usepackage{amsmath,amsthm,amsfonts,amssymb}
\def\qed{{\unskip\nobreak\hfil\penalty50
\hskip2em\hbox{}\nobreak\hfil$\square$
\parfillskip=0pt \finalhyphendemerits=0\par}\medskip}
\def\proof{\trivlist \item[\hskip \labelsep{\bf Proof\ }]}
\def\endproof{\null\hfill\qed\endtrivlist}

\def\isom{\cong}
\def\ti{\tilde}
\def\lan{\langle}
\def\ran{\rangle}

\def\Aut{{\mathrm {Aut}}}

\def\id{{\mathrm {id}}}

\def\Tr{{\mathrm {Tr}}}
\def\Vir{{\mathrm {Vir}}}

\def\a{\alpha}
\def\be{\beta}

\def\epsilon{\varepsilon}

\def\la{\lambda}
\def\La{\Lambda}
\def\phi{\varphi}
\def\si{\sigma}
\def\th{\theta}
\def\om{\omega}

\newtheorem{theorem}{Theorem}[section]
\newtheorem{lemma}[theorem]{Lemma}
\newtheorem{corollary}[theorem]{Corollary}
\newtheorem{definition}[theorem]{Definition}

\newtheorem{proposition}[theorem]{Proposition}
\newtheorem{remark}[theorem]{Remark}
\newtheorem{example}[theorem]{Example}

\def\A{{\cal A}}
\def\B{{\cal B}}
\def\C{{\cal C}}

\def\M{{\cal M}}

\def\Z{{\mathbb Z}}
\def\F{{\mathbb F}}

\title{{\bf Local conformal nets arising from\\
framed vertex operator algebras}}

\author{
{\sc Yasuyuki Kawahigashi}\footnote{Supported in part by JSPS.}\\
Department of Mathematical Sciences\\
University of Tokyo, Komaba, Tokyo, 153-8914, Japan\\
e-mail: {\tt yasuyuki@ms.u-tokyo.ac.jp}\\
\vphantom{X}\\
{\sc Roberto Longo}\footnote{Supported in part by GNAMPA and MIUR.}\\
Dipartimento di Matematica\\
Universit\`a di Roma ``Tor Vergata''\\
Via della Ricerca Scientifica, 1, I-00133 Roma, Italy\\
e-mail: {\tt longo@mat.uniroma2.it}}
\begin{document}
\maketitle

\begin{abstract}
We apply an idea of framed vertex operator algebras to a construction
of local conformal nets of (injective type III$_1$)
factors on the circle corresponding to various lattice vertex operator
algebras and their twisted orbifolds.
In particular, we give a local conformal net corresponding to the
moonshine vertex operator algebras of Frenkel-Lepowsky-Meurman.
Its central charge is 24, it has a trivial representation
theory in the sense that the vacuum sector is the only irreducible
DHR sector, its vacuum character is the modular invariant
$J$-function and its automorphism group (the {\sl gauge group})
is the Monster group.   We use our previous tools such as
$\alpha$-induction and complete rationality to study extensions
of local conformal nets.
\end{abstract}

\section{Introduction}

We have two mathematically rigorous approaches to study chiral
conformal field theory using infinite
dimensional algebraic systems.  One is algebraic quantum
field theory where we study local conformal nets of von Neumann
algebras (factors) on the circle, and the other is theory of vertex
operator algebras.  One local conformal net of factors
corresponds to one vertex operator algebra, at least conceptually,
and each describes one chiral conformal field theory.
Since these two mathematical theories are supposed to
study the same physical objects, it is natural that the two
theories have much in common.  For example, both theories have
mathematical objects corresponding to the affine Lie algebras and
the Virasoro algebra, and also, both have simple current
extension, the coset construction, and the orbifold construction
as constructions of a new object from a given object.
However, the interactions between the two theories have been
relatively small, and different people have studied different
aspects of the two approaches from different motivations.
Comparing the two theories, one easily notices that study of
lattice vertex operator algebras and their twisted orbifolds
has been extensively pursued, but that the corresponding 
study in algebraic quantum field theory is relatively small,
although we have had some works such as
\cite{BMT}, \cite{DX}, \cite{S}, \cite{X3}.  The most
celebrated vertex operator algebra, the moonshine vertex
operator algebra \cite{FLM}, belongs to this class, and almost
nothing has been studied about its counterpart in algebraic quantum
field theory.
The automorphism group of this moonshine vertex operator algebra
is the Monster group \cite{G}, the largest among the 26 sporadic finite
simple groups, and Borcherds has solved the celebrated
moonshine conjecture
of Conway-Norton \cite{CN} on the McKay-Thompson series arising
from this moonshine vertex operator algebra
in \cite{B}.  Our aim in this paper is to study such counterparts
in algebraic quantum field theory,
using an idea of framed vertex operator algebras in
\cite{DGH}.  In particular, we construct a local conformal net
corresponding to the moonshine vertex operator algebra, the
``moonshine net'', in 
Example \ref{moon}.  It has central charge 24 and
a trivial representation theory, its vacuum character is
the modular invariant $J$-function, and its automorphism group is
the Monster group.

Here we briefly explain the relation between local conformal
nets and vertex operator algebras.  In algebraic quantum field
theory, we consider a family of von Neumann algebras of
bounded linear operators generated
by self-adjoint operators corresponding to the ``observables''
in a spacetime region.  (The book \cite{H} is a standard
textbook.)  In the case of chiral conformal field
theory, the ``spacetime'' is the one-dimensional circle and
the spacetime regions are intervals, that is, non-dense non-empty
open connected sets on the circle.  We thus have a family of
von Neumann algebras $\{\A(I)\}_I$ parameterized by intervals
$I$ on a fixed Hilbert space which has a vacuum vector.  Other
properties such as covariance and locality are imposed upon
this family of von Neumann algebras.  Such a family of von Neumann
algebras is called a local conformal net, or simply a net,
of von Neumann algebras, or factors, where a factor means a von
Neumann algebra with a trivial center.  See \cite{GL,KL1} for the
precise set of axioms.  If we start with a vertex operator algebra
(with unitarity), 
then each vertex operator should be an operator-valued distribution
on the circle, so we apply test functions supported on an interval
$I$ and obtain an algebra of bounded operators generated by these
(possibly unbounded) operators arising from such pairing.  This
should give a local conformal nets $\{\A(I)\}_I$, at least
conceptually.  Note that locality for a local conformal net
$\{\A(I)\}_I$ takes a very simple form, that is, if $I$ and
$J$ are disjoint intervals of the circle, then $\A(I)$ and
$\A(J)$ mutually commute.  The axioms of vertex operator algebras
essentially come from Fourier expansions of Wightman fields on
the circle.  By the so-called state-field correspondence, one 
obtains a bijection between a certain set of Wightman fields and 
a dense subspace of the Hilbert space.  With this
identification, each vector in this dense subspace gives a
vertex operator and through the Fourier expansion coefficients,
it produces a countable family of (possibly unbounded)
operators parameterized by integers.
An abstract vertex operator algebras is, roughly speaking, such
a vector space each vector of which gives a countable family
of operators, subject to various compatibility conditions of
these operators such as locality and covariance.  The diffeomorphism
covariance is represented as existence of a special vector called
the Virasoro element. See
\cite{FLM} for the precise definition.  Fredenhagen and J\"or\ss
\cite{FJ} have given a rigorous correspondence between a local
conformal nets of factors and a family of Wightman
fields, but still, the exact relations between local conformal
nets and vertex operator algebras are not clear.  (Note that
unitarity is always assumed from the beginning in algebraic
quantum field theory.)

Dong, Griess and H\"ohn \cite{DGH} have established
a general theory of framed vertex operator algebras and
in particular shown how the moonshine vertex operator algebra
of Frenkel-Lepowsky-Meurman decomposes
into irreducible modules as a module of the 48th
tensor power of the Virasoro vertex operator algebra of
central charge $1/2$. (Also see predating works \cite{D,DMZ}.)
A framed vertex operator algebra is
an extension of a tensor power of the Virasoro vertex operator
algebra having central charge $1/2$.  We ``translate'' this work
into the framework of algebraic quantum field theory in this
paper, which is a nontrivial task.
We hope that the operator algebraic viewpoint
gives a new light on the structure of the Monster group, and
in particular, the Monstrous Moonshine \cite{CN}.

A framed vertex operator algebra is a particular example of 
an extension of a vertex operator algebra.  A general extension
problem of a vertex operator algebra is formulated as follows. 
Take a vertex operator algebra $V$ and let $\{M_i\}_i$ be a set of
representatives from equivalence classes of irreducible
$V$-modules.  We take $M_0$ to be $V$ itself regarded as a
$V$-module.  Then we make a direct sum $\bigoplus_i n_i M_i$,
where $n_i$ is a multiplicity and $n_0=1$, and we would like
to know when this direct sum has a structure of a vertex 
operator algebra with the same Virasoro element in $M_0$
(and how many different vertex operator algebra structures we
have on this).  The corresponding extension problem has been
studied well in algebraic quantum field theory.  That is,
for a local
conformal net $\A$ with irreducible Doplicher-Haag-Roberts
(DHR) sectors $\{\la_i\}_i$ (as in \cite{DHR} and \cite{FRS}),
we would like to determine when $\bigoplus_i n_i \la_i$ gives
an extension of a local conformal net, where $n_i$ is a
multiplicity, $\la_0$ is the vacuum sector, and $n_0=1$.
(A DHR sector of a local conformal net is a unitary equivalence
class of representations of a local conformal net.  A local
conformal net acts on some Hilbert space from the beginning,
and this is called a vacuum representation when regarded as
a representation.  The representation category of a local conformal
net is a braided tensor category.)
A complete characterization when this gives an extension of 
a local conformal net was given in \cite[Theorem 4.9]{LR} using
a notion
of a $Q$-system \cite{L2}.  (An essentially same characterization
was given in \cite[Definition 1.1]{KO} in the framework
of abstract tensor categories.  Conditions 1 and 3 there
correspond to the axioms of the $Q$-system in \cite{L2},
Condition 4 corresponds to irreducibility, and Condition 2
corresponds to chiral locality in \cite[Theorem 4.9]{LR}
in the sense of \cite[page 454]{BEK1}.)
When we do have such an extension,
we can compare representation theories of the original local
conformal net and its extension by a method of $\alpha$-induction,
which gives a (soliton) representation of a larger local conformal
net from a representation of a smaller one, using a braiding.
This induction machinery was defined in \cite{LR}, many
interesting properties and examples, particularly ones related
to conformal embedding, were found in \cite{X1}, and it has
been further studied in \cite{BE,BE4,BEK1,BEK2,BEK3}, particularly
on its relation to modular invariants.  Classification of
all possible structures on a given $\bigoplus_i n_i \la_i$
was studied as a certain 2-cohomology problem in \cite{IK} and
finiteness of the number of possible structures has been proved
in general
there.  Further studies on such 2-cohomology were made in \cite{KL2}
and various vanishing results on this 2-cohomology were obtained.
Using these tools, a complete classification of such extensions
for the Virasoro nets with central charge less than 1, arising
from the usual Virasoro algebra, has been obtained in \cite{KL1}.
The classification list consists of the Virasoro nets themselves,
their simple current extensions of index 2, and four exceptionals
at the central charges $21/22, 25/26, 144/145, 154/155$.
By the above identification of the two extension problems, this
result also gives a complete classification of extensions of
the Virasoro vertex operator algebras with $c<1$.

We refer to \cite{KL1} and its references for a general definition
and properties of local conformal nets (with diffeomorphism
covariance) and $\alpha$-induction.  Also see \cite{K1,K2}
for reviews on \cite{KL1}.  Subfactor theory initiated by Jones
\cite{J} also plays an important role.  See \cite{EK} for 
general theory of subfactors.

\section{Framed vertex operator algebras and extension of local
conformal nets}

Staszkiewicz \cite{S} gave a construction of a local
conformal net from a lattice $L$ and constructed representations
of the net corresponding to the elements of $L^*/L$, where $L^*$ is
the dual lattice of $L$.  Since he also computed the $\mu$-index
of the net, the result in \cite[Proposition 24]{KLM} implies
that the DHR sectors he constructed exhaust all.
Furthermore, after the
first version of this paper was posted on arXiv, a paper of
Dong and Xu \cite{DX} appeared on arXiv studying the operator
algebraic counterparts of lattice vertex operator algebras and
their twisted orbifolds.  Here we use a different approach based
on an idea of framed vertex operator algebras in \cite{DGH}.

We start with a Virasoro net with $c=1/2$ constructed as
a coset in \cite{X2} and studied in \cite{KL1}.  In
\cite{KL1}, we have shown that this net is completely rational 
in the sense of \cite[Definition 8]{KLM}, and that it
has three irreducible DHR sectors of statistical
dimensions $1, \sqrt2, 1$, respectively.  So the
$\mu$-index, the index of the two-interval subfactor as in
\cite[Proposition 5]{KLM}, of this net is $4$ and their
conformal weights are $0,1/16,1/2$, respectively.
(The $\mu$-index of a local conformal net is equal to
the square sum of the statistical dimensions of all irreducible
DHR sectors of the net by \cite[Proposition 24]{KLM}.  This quantity
of a braided tensor category also plays an important role in
studies of quantum invariants in three-dimensional topology.)
We denote this net by $\Vir_{1/2}$.  We next consider the net
$\Vir_{1/2}\otimes\Vir_{1/2}$.  This net is also completely
rational and has $\mu$-index $16$.  It
has $9$ irreducible DHR sectors, since each such sector
is a tensor product of two irreducible DHR sectors of
$\Vir_{1/2}$.  Using the conformal weights of the two such
irreducible DHR sectors of $\Vir_{1/2}$, we label the $9$
irreducible DHR sectors of $\Vir_{1/2}\otimes\Vir_{1/2}$ as
$$\la_{0,0}, \la_{0,1/16}, \la_{0,1/2},
\la_{1/16,0}, \la_{1/16,1/16}, \la_{1/16,1/2},
\la_{1/2,0}, \la_{1/2,1/16}, \la_{1/2,1/2}.$$

We denote the conformal weight, the lowest eigenvalue of the 
image of $L_0$,
and the conformal spin of an irreducible sector $\la$ by $h_\la$
and $\om_\la$, respectively.  By the spin-statistics theorem of
Guido-Longo \cite{GL}, we have $\exp(2\pi i h_\la)=\om_\la$.

\begin{lemma}
\label{spin1}
Let $\A$ be a local conformal net on $S^1$.  Suppose we have
a finite system $\{\la_j\}_j$ of irreducible DHR endomorphisms of
$\A$ and each $\la_j$ has a statistical dimension $1$.
If all $\la_j$ have a conformal spin $1$, then the crossed product
of $\A$ by the finite abelian group $G$ given by $\{\la_j\}_j$
produces a local extension of the net $\A$.
\end{lemma}

\begin{proof}
First by \cite[Lemma 4.4]{R2}, we can change the representatives
$\{\la_j\}_j$ within their unitary equivalence classes so that
the system $\{\la_j\}_j$ gives an action of $G$ on a local algebra
$A(I)$ for a fixed interval $I$.  That is, we decompose $G=\prod_i G_i$,
where each $G_i$ is a cyclic group, and localize generators of
$G_i$'s on mutually disjoint intervals within $I$.  By the 
assumption on the conformal spins, we can adjust each generator
so that it gives an action of $G_i$.  Then as in
\cite[Part II, Section 3]{BE}, we can make an extension of the
local conformal net $\A$ to the crossed product by the $G$-action.
The proof of \cite[Part II, Lemma 3.6]{BE} gives the desired
locality, because the conformal spins are now all 1.
(Lemma 3.6 in \cite[Part II]{BE} deals with only actions of a
cyclic group, but the same argument works in our current setting.)
\end{proof}

This extended local conformal net is called a {\sl simple current
extension} of $\A$ by $G$.
We now apply Lemma \ref{spin1} to the net
$\Vir_{1/2}\otimes\Vir_{1/2}$ and its system of irreducible
DHR sectors consisting of $\la_{0,0}, \la_{1/2,1/2}$.
Note that the conformal weight of the sector $\la_{1/2,1/2}$ is $1$,
so its conformal spin is $1$ and we can apply Lemma \ref{spin1}
with $G=\Z_2$. 

\begin{proposition}\label{U1-4}
This net $\A$ is completely rational  and it 
has four irreducible DHR sectors and all have
statistical dimensions $1$.  The fusion rules of these four sectors are
given by the group $\Z_4$ and the conformal
weights are $0,1/16,1/2,1/16$.
\end{proposition}

\begin{proof}
The net $\A$ is completely rational by \cite{L3}.
By \cite[Proposition 24]{KLM}, the $\mu$-index is $4$.
Since $\a^\pm$-induction of $\la_{0,0}$ and $\la_{0,1/2}$ give
two irreducible DHR sectors of the net $\A$, we obtain
two inequivalent irreducible DHR sectors of $\A$ having statistical
dimensions $1$.  The $\a^\pm$-induction of $\la_{1/16,1/16}$ 
gives either one irreducible DHR sector of statistical dimension
$2$ or two inequivalent  irreducible DHR sectors of statistical 
dimensions $1$.  Since the square sum of the statistical dimensions
of all the irreducible DHR sectors is $4$, the latter case in the
above occurs.  In this way, we obtain four inequivalent irreducible
DHR sectors of statistical dimensions $1$, and these exhaust all
the irreducible DHR sectors of $\A$.  Note that two of
the irreducible DHR sectors have conformal weights $1/8$, and thus
conformal spins $\exp(\pi i/4)$.  

The fusion rules are
given by either $\Z_4$ or $\Z_2\times\Z_2$.
If we have the fusion rules of $\Z_2\times\Z_2$,
then the two irreducible DHR sectors of conformal weights $1/8$ must
have order 2, but this would violate
\cite[Corollary on page 343]{R1}.  Thus we obtain the former
fusion rules.
\end{proof}

\begin{remark}{\rm
This net $\A$ has a central charge $1$.  In \cite{X5}, Xu classified
all local conformal nets on $S^1$ under an additional assumption called
a spectrum condition.  It is conjectured that this spectrum condition
always holds.  Thus, our net $\A$ should be in the classification
list in \cite{X5}.  Since we have only one net having the $\mu$-indexes
equal to $1$ in the list in \cite{X5}, our net $\A$ should be
isomorphic to the net $U(1)_4$.  (See \cite[Section 3.5]{X4} for
the nets $U(1)_{2k}$.)
}\end{remark}

We now use the framework of Dong-Griess-H\"ohn \cite{DGH}.
As \cite[Section 4]{DGH}, let $C$ be a doubly-even linear
binary code of length $d\in8\Z$ containing the vector
$(1,1,\dots,1)$.  (A doubly-even code is sometimes called
type II.)  This $C$ is a subset of $\F_2^d$, where
$\F_2$ is a field consisting two elements, and any
element in $C$ is naturally regarded as an element in $\Z^d$.
As in \cite{CS}, we associate two even lattices with such
a code $C$ as follows.
\begin{eqnarray*}
L_C&=&\{(c+x)/\sqrt2\mid c\in C, x\in(2\Z)^d\},\\
\ti L_C&=&\{(c+y)/\sqrt2\mid c\in C, y\in (2\Z)^d, 
\sum y_i\in 4\Z\}\cup\\
&&\{(c+y+(1/2,1/2,\dots,1/2))/\sqrt2\mid
c\in C, y\in (2\Z)^d, 1-(-1)^{d/8}+\sum y_i\in 4\Z\}.
\end{eqnarray*}
We then have the corresponding vertex operator algebras
$V_{L_C}, V_{\ti L_C}$ as in \cite{FLM}, but the construction
directly involves vertex operators and it seems very difficult
to ``translate'' this construction into the operator 
algebraic framework.  So we will take a different approach
based on an idea of a framed vertex operator algebras.
(If the lattice is $D_1$, we do have a
counterpart of the lattice vertex operator algebra
$V_{D_1}$ and it is the above net $\A$.  This has been
already noted in \cite[page 14072]{X3} and also is a basis of
\cite{S}.)

Let $L$ be a positive definite even
lattice of rank $d$ containing $D_1^d$ as a sublattice, where
$D_1=2\Z$ and $\lan \a,\be\ran=\a\be$ for $\a,\be\in2\Z$.  
Such a sublattice is called a $D_1$-frame.
Recall that we have a non-degenerate symmetric $\Z$-bilinear
form $\lan\cdot,\cdot\ran$ on $L$ and $\lan \a, \a\ran\in2\Z$
for all $\a\in L$.  Note that $(D_1^*/D_1)^d$ is isomorphic
to $\Z_4^d$, thus we have a set $\Delta(L)=L/D_1^d\subset
(D_1^*/D_1)^d$ and this set is a code over $\Z_4$, that is, 
$\Delta(L)$ is a subgroup of $\Z_4^d$.  An element of a code
is called a codeword.  It is known that
$\Delta(L)$ is self-dual if and only if $L$ is
self-dual.  (See \cite[page 426]{DGH}.)
For the above lattices $L_C, \ti L_C$, the
corresponding $\Z_4$-codes are explicitly known as follows.
(Also see \cite[page 426]{DGH}.)

Let $\hat{}$ be a map from $\F_2^2$ to $\Z_4^2$ given by
$00\mapsto 00$, $11\mapsto 20$, $10\mapsto 11$, and
$01\mapsto 31$.  By regarding $\F_2^d=(\F_2^2)^{d/2}$ and
$\Z_4^d=(\Z_4^2)^{d/2}$, and applying the map $\hat{}$
componentwise, we obtain the map, still denoted by $\hat{}$,
from $\F_2^d$ to $\Z_4^d$.  We also define
$\Sigma_2^n=\{(00),(22)\}^n$ and let $(\Sigma_2^n)_0$ be the
subcode of the $\Z_4$-code $\Sigma_2^n$ consisting of
codewords of Hamming weights divisible by 4.  (The Hamming
weight of a codeword is the number of nonzero entries of
the codeword.)  The we have
\begin{eqnarray*}
\Delta(L_C)&=&\hat C+\Sigma_2^{d/2},\\
\Delta(\ti L_C)&=& \hat C +(\Sigma_2^{d/2})_0 \cup
\hat C +(\Sigma_2^{d/2})_0+\left\{
\begin{array}{ll}
(1,0,\cdots,1,0,1,0),&{\rm if\ }d\equiv 0 \mod 16,\\
(1,0,\cdots,1,0,3,2),&{\rm if\ }d\equiv 8 \mod 16.
\end{array}\right.
\end{eqnarray*}

Let $G$ be the abelian group given by $\Delta(L_C)$ or
$\Delta(\ti L_C)$, regarded as a subgroup of $\Z_4^d$.
We consider a completely rational local conformal net
$\A^{\otimes d}$.  As above, each irreducible DHR
sector of this net is labeled with an element of $\Z_4^d$
naturally.

\begin{lemma}\label{spin2}
Each irreducible DHR sector corresponding to an element
in $G$ has a conformal spin $1$.
\end{lemma}

\begin{proof}
We consider a conformal weight of each DHR sector.  For the four
elements
$00, 01, 10, 11 \in \F_2^2$, the images by the map $\hat{}$
are $00, 31, 11, 20$, respectively.  For a net $\A\otimes\A$, 
we have the corresponding irreducible DHR sectors,
and their conformal weights are $0,1/4,1/4,1/2$, respectively.  This
means that $1\in \F_2$ has a contribution $1/4$ to the
conformal weight and $0\in F_2$ has a contribution $0$.  Since the code
$C$ is doubly-even, the number of nonzero entries in
any element $c$ in $C$  is a multiple of $4$.  This means that
for any element $c\in C$, the irreducible DHR sector
corresponding to $\hat c$ has an integer conformal weight.

Since the irreducible DHR sector corresponding to the
element $2\in \Z_4$ has a conformal weight $1/2$, the
irreducible DHR sector corresponding to  any element in
$\Sigma_2^{d/2}$ also has an integer conformal weight.
The irreducible DHR sector corresponding to the
elements $1, 3\in \Z_4$ have conformal weights $1/8$.  So,
if $d\equiv0\mod 16$, then  the irreducible DHR sector
corresponding to the element $(1,0,\cdots,1,0,1,0)\in \Z_4^d$
has a conformal weight $d/16$, which is an integer.
If $d\equiv8\mod 16$, then  the irreducible DHR sector
corresponding to the element $(1,0,\cdots,1,0,3,2)\in \Z_4^d$
has a conformal weight $(d+8)/16$, which is again an integer.

We have thus proved that for any element in $\Delta(L_C)$ or
$\Delta(\ti L_C)$, the corresponding irreducible DHR sector
has an integer conformal weight, hence a conformal spin $1$ by the
spin-statistics theorem \cite{GL}.
\end{proof}

We then have the following theorem.

\begin{theorem}\label{holo}
We can extend $\A^{\otimes d}$ to a local conformal net
as a simple current extension by the group $G=\Delta(L_C),
\Delta(\ti L_C)$ as above.  If the code $C$ is self-dual,
the extended local conformal net has $\mu$-index $1$.
\end{theorem}

\begin{proof}
The first claim is clear by Lemmas \ref{spin1}, \ref{spin2}.
If the code $C$ is self-dual, its cardinality is $2^{d/2}$,
so the cardinality of $G$ is $2^d$.  (See \cite[(10.1.8)]{FLM}.)
Since the $\mu$-index
of $\A^{\otimes d}$ is $4^d$, we obtain the conclusion
by \cite[Proposition 24]{KLM}.
\end{proof}

The condition that the $\mu$-index is 1 means that the 
vacuum representation is the only
irreducible representation of the net.  This property
is well-studied in theory of vertex operator algebras and
such a vertex operator algebra is called holomorphic (or also
self-dual).  So
we give the following definition.

\begin{definition}{\rm
A local conformal net is said to be {\sl holomorphic} if
it is completely rational and has $\mu$-index $1$.
}\end{definition}

Recall that the vacuum character of a local conformal net is
defined as $$\Tr(e^{2\pi i\tau (L_0-c/24)})$$
on the vacuum Hilbert space.
The vacuum character of a vertex operator algebra has been
defined in a similar way.

\begin{proposition}\label{corresp}
Let $L$ be an even lattice arising from a doubly-even linear
binary code $C$ of length $d\in8\Z$ containing the vector
$(1,1,\dots,1)$ as $L_C$ or $\ti L_C$ as above.  Denote the
local conformal net arising from $L$ as above by $\A_L$
and the vertex operator algebra arising from $L$
by $V_L$.  If we start with a self-dual code $C$,
both $\A_L$  and $V_L$ are holomorphic.
The central charges of $\A_L$ and
$V_L$ are both equal to $d$.  Furthermore, the vacuum characters of
$\A_L$ and $V_L$ are equal.
\end{proposition}

\begin{proof}
The statements on holomorphy and the central charge are now obvious.
The Virasoro net $\Vir_{1/2}$ and the Virasoro vertex operator
algebra $L(1/2,0)$ have the same representation theory 
consisting of three irreducible objects, and
their corresponding characters are also equal.
So both of their $2d$-th tensor powers also have the 
same representation theory and the same characters.
When we decompose $V_L$ as a module over $L(1/2,0)^{\otimes 2d}$
and decompose the canonical endomorphism of the inclusion
of $\Vir_{1/2}^{\otimes 2d}\subset \A_L$   into irreducible
DHR sectors of $\Vir_{1/2}^{\otimes 2d}$, we can identify the
two decompositions from the above construction.  This implies,
in particular, that the vacuum characters are equal.
\end{proof}

We study the characters for local conformal net $\A$.
For an irreducible DHR sector $\la$ of $\A$, we define the
specialized character
$\chi_\lambda(\tau)$ for complex numbers $\tau$ with
${\rm Im}\; \tau > 0$ as follows.
$$\chi_\lambda(\tau)=\Tr(e^{2\pi i\tau (L_{0,\la}-c/24)}),$$
where the operator $L_{0,\la}$ is conformal Hamiltonian in 
the representation $\la$ and $c$ is the central charge.
For a module $\la$ of vertex operator algebra $V$, we define
the specialized character by the same formula.  Zhu \cite{Z}
has proved under some general condition that the group
$SL(2,\Z)$ acts on the linear span of these specialized
characters through the usual change of variables $\tau$
and this result applies to lattice vertex operator algebras.
In particular, if a lattice vertex operator algebra is
holomorphic, then the only specialized character is the
vacuum one, and this vacuum character must be invariant under
$SL(2,\Z)$, that is, a modular function.  We thus know that
our local conformal nets $\A_L$ also have modular functions
as the vacuum characters.  This can be proved also as follows.
For the Virasoro net $\Vir_{1/2}$, the specialized characters
for the irreducible DHR sectors coincide with the usual characters
for the modules of the Virasoro vertex operator algebra $L(1/2,0)$.
Then by Proposition 2 in \cite{KL3}, the vacuum character of
a holomorphic net $\A_L$ is a modular invariant function.

We show some concrete examples, following \cite[Section 5]{DGH}.

\begin{example}{\rm
The Hamming code $H_8$ is a self-dual double-even binary
code of length $8$.  This gives the $E_8$ lattice as both
$L_{H_8}$ and $\ti L_{H_8}$, and the above construction then 
produces a holomorphic local conformal net
with central charge $8$.  In theory of vertex operator algebras,
the corresponding lattice vertex operator algebra $V_{E_8}$
is isomorphic to the one arising from the affine Lie
algebra $E^{(1)}_8$ at level $1$, thus we also expect that
the above local conformal net is isomorphic to the
one arising from the loop group representation for $E^{(1)}_8$
at level $1$.
}\end{example}

\begin{example}\label{golay}{\rm
Consider the Golay Code ${\cal G}_{24}$ as in
\cite[Chapter 10]{FLM}, which is a
self-dual double-even binary code of length $24$.
Then the lattice $\ti L_{{\cal G}_{24}}$ is the Leech lattice,
which is the unique positive definite even self-dual
lattice of rank 24 having no vectors of square length $2$
as in \cite{FLM}.
Its vacuum character is the following, which is equal to
$J(q)+24$.
$$q^{-1}+24+196884q+21493760q^2+8642909970 q^3+\cdots.$$
}\end{example}

\section{Twisted orbifolds and $\alpha$-induction}

We now perform the twisted orbifold construction for
the net arising from a lattice as above.
Let $C$ be self-dual double-even binary code of length $d$
and $L$ be $L_C$ or $\ti L_C$, the lattice arising from  $C$ as
in the above section.  We study the net $\A_L$, which is holomorphic
by Theorem \ref{holo}.

The basic idea of the twisted orbifold construction can be
formulated in the framework of algebraic quantum field theory
quite easily.  We follow Huang's idea \cite{Hu} in theory of
vertex operator algebras.  Take a holomorphic local conformal
net $\A$ and its automorphism $\si$ of order 2.  Then the
fixed point net $\A^\si$ has $\mu$-index 4, and it is  easy
to see that this fixed point net has four irreducible DHR
sectors of statistical dimensions 1, and the fusion rules
are given by $\Z_2\times\Z_2$.  Besides the vacuum sector,
two of them have conformal spins 1, so we can make a simple
current extension of $\A^\si$ with each of these two sectors.
One of them obviously gives back $\A$, and the other is
the twisted orbifold of $\A$.  (Here we have presented only
an outline.  The actual details are given below.)  In this
way, we can obtain the counterpart of the moonshine vertex
operator algebra from $\A_\La$, where $\La$ is the Leech
lattice, but it is hard to see the structure of the twisted
orbifold in this abstract approach.  For example, an important
property of the moonshine vertex operator algebras is that its
vacuum character is the modular invariant $J$-function, but
it is hard to see this property for the corresponding net in
the above approach.  Moreover, we are interested in the
automorphism group of the twisted orbifold, but it is
essentially impossible to determine the automorphism group
through this abstract construction of a local conformal net.
So we need a more concrete study of the twisted
orbifold construction for local conformal nets, and we will
make such a study with technique of $\a$-induction below.

Recall that the net $\A$ is constructed as a crossed
product by a $\Z_2$-action in Proposition \ref{U1-4}, 
so it has a natural dual
action of $\Z_2$.  We take its $d$-th tensor power on
${\A}^{\otimes d}$ and denote it by $\sigma$.  This is an
automorphism of the net ${\A}^{\otimes d}$ of order 2.
We extend ${\A}^{\otimes d}$ to a local conformal net $\A_L$ as in
Theorem \ref{holo}.   Then $\sigma$ extends to an
automorphism of $\A_L$ of order 2 by Lemma \ref{spin2} and
\cite[Proposition 2.1]{BDLR}, and we denote the extension 
again by $\sigma$.  Note that the fusion rules of
${\A}^{\otimes d}$ are naturally given by $\Z_4^d$, and the
action of $\si$ on the fusion rules $\Z_4$ is given by
$j\mapsto -j\in\Z_4$.

\begin{lemma}\label{L1}
The fixed point net $\A_L^\si$ has four
irreducible DHR sectors of statistical dimensions $1$ and
they exhaust all the irreducible DHR sectors.
\end{lemma}

\begin{proof}
The $\mu$-index of the fixed point net is $4$ by 
\cite[Proposition 24]{KLM}.  So either we have four
irreducible DHR sectors of statistical dimensions $1$
or we have three irreducible DHR sectors of statistical
dimensions $1$, $1$, $\sqrt2$.  Suppose we had the
latter case and draw the induction-restriction graph
for the inclusion $\A_L^\si\subset \A_L$.  Then the index 
of this inclusion
is $2$ and the Dynkin diagram $A_3$ is the only graph
having the Perron-Frobenius eigenvalue $\sqrt2$,  and
statistical dimensions are always larger than or equal
to $1$, thus we would have an irreducible DHR sector
of statistical dimension $\sqrt2$ for the holomorphic
net $\A_L$, which is impossible. We thus know that
we have four irreducible DHR sectors of statistical 
dimensions $1$ and these give all the irreducible DHR sectors
of the net $\A_L^\sigma$.
\end{proof}

We study the following inclusions of nets.
$$\begin{array}{ccccc}
&&\A^{\otimes d} &\subset & \A_L\\
&&\cup && \cup \\
\Vir_{1/2}^{\otimes 2d}&\subset&
(\A^{\otimes d})^\si &\subset & \A_L^\si.
\end{array}$$

Our aim is to study decompositions of irreducible DHR sectors
of the net $\A_L^\si$ restricted to $\Vir_{1/2}^{\otimes 2d}$.
We describe such decompositions through studies of soliton
sectors of $\A^{\otimes d}$ and $\A_L$.  This is because 
the inclusions $\Vir_{1/2}^{\otimes 2d}\subset \A^{\otimes d} 
\subset \A_L$ are easier to study.

\begin{lemma}\label{L2}
The irreducible DHR sectors of the net $(\A^{\otimes d})^\si$
consist of $4^{d-1}$ sectors of statistical dimensions $2$,
$4^d$ sectors of statistical dimensions $1$, 
and $2^{d+1}$ sectors of statistical dimensions $2^{d/2}$.

After applying the $\a^+$-induction for the inclusion
$(\A^{\otimes d})^\si\subset \A^{\otimes d}$ and the irreducible
DHR sectors of $(\A^{\otimes d})^\si$ as above,
we obtain irreducible soliton sectors consisting of
$4^d$ sectors of statistical dimensions $1$ and
$2^d$ sectors of statistical dimensions $2^{d/2}$.
The $4^d$ sectors of statistical dimensions $1$ precisely give
the irreducible DHR sectors of the net $\A^{\otimes d}$.
\end{lemma}

\begin{proof}
Let
$$K=\{(a_j)\in \Z_2^d\mid \sum_j a_j=0\in\Z_2\},$$
which is a subgroup of $\Z_2^d$ of order $2^{d-1}$.
Then the inclusion of nets $\Vir_{1/2}^{\otimes 2d}\subset
(\A^{\otimes d})^\si\subset \A^{\otimes d}$ is identified
with $\Vir_{1/2}^{\otimes 2d}\subset
(\Vir_{1/2}^{\otimes 2d})\rtimes K\subset
\Vir_{1/2}^{\otimes 2d}\rtimes \Z_2^d$.  From this
description of the inclusion and information of
the representation category of $\Vir_{1/2}^{\otimes 2d}$,
we obtain the description of the irreducible DHR sectors
for the net $(\A^{\otimes d})^\si$ as follows.

Take an irreducible DHR sector $\la=\bigotimes_{k=1}^d
\la_k$ of $\Vir_{1/2}^{\otimes 2d}$, where each $\la_k$ is
one of the irreducible DHR sectors
$\la_{0,1/16}, \la_{1/2,1/16}, \la_{1/16,0}, \la_{1/16,1/2}$
of $\Vir_{1/2}\otimes \Vir_{1/2}$.  Let $\th$ be the dual 
canonical endomorphism for the extension 
$\Vir_{1/2}^{\otimes 2d}\subset
(\Vir_{1/2}^{\otimes 2d})\rtimes K$.  Then we see that the
monodromy of $\th$ and $\la$ is trivial because each automorphism
appearing in $\th$ has a trivial monodromy with $\la$ by the
description of the $S$-matrix given by \cite[(10.138)]{DMS}.
Then \cite[Part I, Proposition 3.23]{BE} implies that $\a^+$-
and $\a^-$-inductions of this $\la$ give the same sector, thus
they produce a DHR sector of the net
$(\Vir_{1/2}^{\otimes 2d})\rtimes K$.  It is irreducible 
by \cite[Theorem 3.3]{X1} and has a statistical dimension
$2^{d/2}$.  Again by \cite[Theorem 3.3]{X1}, we conclude that
we have $2^{d+1}$ mutually distinct irreducible DHR sectors of 
this form.

We now determine the induction-restriction graph for
the inclusion $(\A^{\otimes d})^\si \subset \A^{\otimes d}$.
By chiral locality, each connected component of this graph for
all the irreducible DHR sectors of $(\A^{\otimes d})^\si$
must be one of the two components in Figure \ref{induc}.

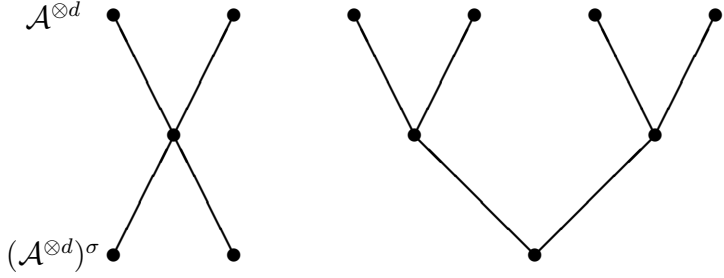
\begin{figure}[tb]
\unitlength 0.8mm
\thicklines
\begin{center}
\begin{picture}(120,60)
\put(10,10){\line(1,2){10}}
\put(30,10){\line(-1,2){10}}
\put(20,30){\line(-1,2){10}}
\put(20,30){\line(1,2){10}}
\put(80,10){\line(-1,1){20}}
\put(80,10){\line(1,1){20}}
\put(60,30){\line(-1,2){10}}
\put(60,30){\line(1,2){10}}
\put(100,30){\line(-1,2){10}}
\put(100,30){\line(1,2){10}}
\put(10,10){\circle*{2}}
\put(30,10){\circle*{2}}
\put(80,10){\circle*{2}}
\put(20,30){\circle*{2}}
\put(60,30){\circle*{2}}
\put(100,30){\circle*{2}}
\put(10,50){\circle*{2}}
\put(30,50){\circle*{2}}
\put(50,50){\circle*{2}}
\put(70,50){\circle*{2}}
\put(90,50){\circle*{2}}
\put(110,50){\circle*{2}}
\put(0,10){\makebox(0,0){$(\A^{\otimes d})^\sigma$}}
\put(0,50){\makebox(0,0){$\A^{\otimes d}$}}
\end{picture}
\end{center}
\caption{Induction-restriction graphs}
\label{induc}
\end{figure}

Since the irreducible DHR sectors of the net $\A^{\otimes d}$
is labeled with elements of $\Z_4^d$ and the action of
$\sigma$ on this fusion rule algebra is given by $-1$ on
each $\Z_4$, the induction-restriction graph for the
irreducible DHR sectors of the net  $\A^{\otimes d}$ contains
$4^{d-1}$ copies of the graphs in Figure \ref{induc1}.
The labels $\a^\pm$ mean that these irreducible DHR sectors
arise from both $\a^+$- and $\a^-$-inductions, and the numbers
represent the statistical dimensions.

\begin{figure}[tb]
\unitlength 0.9mm
\thicklines
\begin{center}
\begin{picture}(160,60)
\put(40,10){\line(-1,1){20}}
\put(40,10){\line(1,1){20}}
\put(90,10){\line(1,2){10}}
\put(110,10){\line(-1,2){10}}
\put(130,10){\line(1,2){10}}
\put(150,10){\line(-1,2){10}}
\put(20,30){\line(-1,2){10}}
\put(20,30){\line(1,2){10}}
\put(60,30){\line(-1,2){10}}
\put(60,30){\line(1,2){10}}
\put(100,30){\line(-1,2){10}}
\put(100,30){\line(1,2){10}}
\put(140,30){\line(-1,2){10}}
\put(140,30){\line(1,2){10}}
\put(40,10){\circle*{2}}
\put(90,10){\circle*{2}}
\put(110,10){\circle*{2}}
\put(130,10){\circle*{2}}
\put(150,10){\circle*{2}}
\put(20,30){\circle*{2}}
\put(60,30){\circle*{2}}
\put(100,30){\circle*{2}}
\put(140,30){\circle*{2}}
\put(10,50){\circle*{2}}
\put(30,50){\circle*{2}}
\put(50,50){\circle*{2}}
\put(70,50){\circle*{2}}
\put(90,50){\circle*{2}}
\put(110,50){\circle*{2}}
\put(130,50){\circle*{2}}
\put(150,50){\circle*{2}}
\put(40,5){\makebox(0,0){$2$}}
\put(90,5){\makebox(0,0){$1$}}
\put(110,5){\makebox(0,0){$1$}}
\put(130,5){\makebox(0,0){$1$}}
\put(150,5){\makebox(0,0){$1$}}
\put(10,55){\makebox(0,0){$1$}}
\put(50,55){\makebox(0,0){$1$}}
\put(90,55){\makebox(0,0){$1$}}
\put(130,55){\makebox(0,0){$1$}}
\put(10,60){\makebox(0,0){$\alpha^\pm$}}
\put(50,60){\makebox(0,0){$\alpha^\pm$}}
\put(90,60){\makebox(0,0){$\alpha^\pm$}}
\put(130,60){\makebox(0,0){$\alpha^\pm$}}
\put(2,10){\makebox(0,0){$(\A^{\otimes d})^\sigma$}}
\put(2,50){\makebox(0,0){$\A^{\otimes d}$}}
\end{picture}
\end{center}
\caption{DHR sectors and $\alpha$-induction (1)}
\label{induc1}
\end{figure}
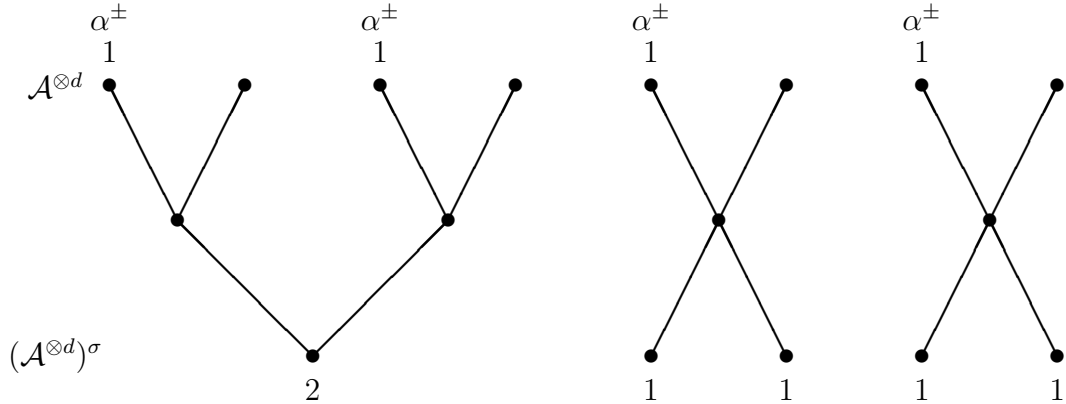

By the above description of the irreducible DHR sectors of 
the net $(\A^{\otimes d})^\si$, we know
that the induction-restriction graph for the inclusion
$(\A^{\otimes d})^\si\subset \A^{\otimes d}$
and all these irreducible DHR sectors of the net $(\A^{\otimes d})^\si$
contains $2^d$ copies of the graph in Figure \ref{induc2},
besides $4^{d-1}$ copies of the graphs in Figure \ref{induc1}
described above.  The labels $\a^+$ and $\a^-$ mean that these
two irreducible DHR sectors result from $\a^+$- and $\a^-$-inductions,
respectively.  (Note that
$4^{d-1}$ copies of the graphs in Figure \ref{induc1} already
give all the irreducible DHR sectors of the net  $(\A^{\otimes d})$,
thus $\a^+$- and $\a^-$-inductions in Figure \ref{induc2} must
produce distinct soliton sectors.)

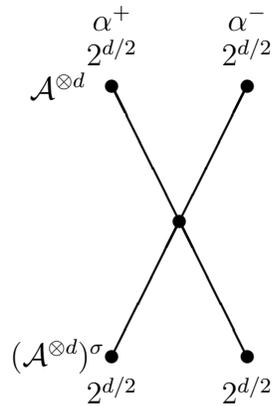
\begin{figure}[tb]
\unitlength 0.9mm
\thicklines
\begin{center}
\begin{picture}(40,60)
\put(10,10){\line(1,2){10}}
\put(30,10){\line(-1,2){10}}
\put(20,30){\line(-1,2){10}}
\put(20,30){\line(1,2){10}}
\put(10,10){\circle*{2}}
\put(30,10){\circle*{2}}
\put(20,30){\circle*{2}}
\put(10,50){\circle*{2}}
\put(30,50){\circle*{2}}
\put(10,5){\makebox(0,0){$2^{d/2}$}}
\put(30,5){\makebox(0,0){$2^{d/2}$}}
\put(10,55){\makebox(0,0){$2^{d/2}$}}
\put(30,55){\makebox(0,0){$2^{d/2}$}}
\put(10,60){\makebox(0,0){$\alpha^+$}}
\put(30,60){\makebox(0,0){$\alpha^-$}}
\put(2,10){\makebox(0,0){$(\A^{\otimes d})^\sigma$}}
\put(2,50){\makebox(0,0){$\A^{\otimes d}$}}
\end{picture}
\end{center}
\caption{DHR sectors and $\alpha$-induction (2)}
\label{induc2}
\end{figure}

The description of these graphs now gives the desired conclusion
of the Lemma, because the total $\mu$-index of the sectors already
described gives the correct $\mu$-index.
\end{proof}

We fix an interval $I\subset S^1$ and consider the following
commuting square of factors.
$$\begin{array}{ccc}
\A^{\otimes d}(I) &\subset & \A_L(I)\\
\cup && \cup \\
(\A^{\otimes d})^\si(I) &\subset & \A_L^\si(I).
\end{array}$$
As $(\A^{\otimes d})^\si(I)$-$(\A^{\otimes d})^\si(I)$
sectors, we consider those arising from the irreducible DHR
sectors of the net $(\A^{\otimes d})^\si$ and draw the
induction-restriction graphs for the above commuting square.

Two of the four irreducible DHR sectors of the net $\A_L^\si$
appear in the decomposition of the vacuum sector of $\A_L$.
Denote the other two irreducible DHR sectors of the net $\A_L^\si$
by $\be_1, \be_2$.  Apply the $\a^+$-induction for the inclusion
$\A_L^\si\subset \A_L$ then both $\be_1, \be_2$ give the same
irreducible soliton sector of statistical dimension $1$.  
Let $\ti\be_1$ be this soliton sector.  We consider the decomposition
of $\ti\be_1$ restricted to $\A^{\otimes d}$.  This gives a direct
sum of irreducible soliton sectors of $\A^{\otimes d}$ arising
from the $\a^+$-induction for the inclusion $(\A^{\otimes d})^\si
\subset \A^{\otimes d}$.  Note that the automorphism group $G$
as in Lemma \ref{spin1} gives a fusion rule subalgebra of the
irreducible DHR sectors of the net $\A^{\otimes d}$.  This group
acts on the system of irreducible soliton sectors arising from 
$\a^+$-induction for the inclusion $(\A^{\otimes d})^\si
\subset \A^{\otimes d}$ by multiplication.  Then this system
of such soliton sectors decomposes into $G$-orbits.

\begin{lemma}\label{L4}
The decomposition of $\ti\be_1$ into such soliton sectors
give one $G$-orbit with a common multiplicity for all the soliton
sectors.
\end{lemma}

\begin{proof}
The subnet $(\A^{\otimes d})^\si\subset \A_L$ is given by 
a fixed point net of a finite group, so all the soliton sectors
of $\A_L$ arising from $\a^\pm$-inductions have statistical
dimension 1 by \cite[Corollary 3.20]{Mu}.
Consider the inclusion $Q/D_1^d\subset L/D_1^d=G\subset \Z_4^d$
as in \cite[page 428]{DGH}.  Let $H$ be the kernel of the action
of $G$ on the set of above irreducible soliton sectors of
statistical dimensions $2^{d/2}$.  We then have
$Q/D_1^d\subset H \subset G$.  Then each connected component
of the induction-restriction graph for the inclusion
$\A^{\otimes d}\subset \A_L$ produces a $G$-orbit of 
such soliton sectors because the dual canonical endomorphism for
the inclusion $\A^{\otimes d}\subset \A_L$ gives a group $G$.
Since all the soliton sectors of $\A_L$ now have dimensions 1,
a simple Perron-Frobenius argument gives that one connected
component of the induction-restriction graph involving the
irreducible soliton sectors of statistical dimensions $2^{d/2}$
of $\A^{\otimes d}$ has one $G$-orbit.  Then a decomposition of
$\ti\be_1$ into soliton sectors
gives such a $G$-orbit with a common multiplicity $[H:Q/D_1^d]$.
\end{proof}

\begin{lemma}\label{L5}
We decompose $\be_1\oplus\be_2$ for the restriction of
$\A_L^\si$ to $\Vir_{1/2}^{\otimes 2d}$.  Then this
decomposition exactly corresponds to the decomposition of
$V_L^T$ in \cite[Proposition 4.8]{DGH}.
\end{lemma}

\begin{proof}
Recall that the irreducible DHR sectors of the net 
$\Vir_{1/2}^{\otimes 2d}$ are labeled
with $J\in\{0,1/2,1/16\}^{2d}$.  So the decomposition of 
$\be_1\oplus\be_2$ for the restriction to
$\Vir_{1/2}^{\otimes 2d}$ is given as $\bigoplus_J n_J \la_J$,
where $n_J$ is the multiplicity of $\la_J$ given by the label
$J$ as above.  The irreducible modules of the vertex operator
algebra $L(1/2,0)^{\otimes 2d}$ are also labeled with the same
$J$ and we denote these irreducible modules as $M_J$.  Then
$V_L^T$ is decomposed as $\bigoplus_J m_J M_J$ as
a $L(1/2,0)^{\otimes 2d}$-module as in \cite[Proposition 4.8]{DGH}.
By comparing this decomposition and Lemma \ref{L4},
we conclude that $n_J=m_J$.
\end{proof}

\begin{lemma}\label{L6}
Among the two irreducible DHR sectors $\be_1,\be_2$
of the net $\A_L^\si$, one has a conformal
spin $1$ and the other has $-1$.  The fusion rules
of all the four irreducible DHR sectors of the net $(\A_L^\si)$
are given by the group $\Z_2\times \Z_2$.
\end{lemma}

\begin{proof}
The weights of $\be_1\oplus\be_2$ consist of integers
and half-integers by Lemma \ref{L5} and the description 
in \cite[Proposition 4.8]{DGH}, since weights are the same
in both setting of the  Virasoro algebra as local conformal nets 
and vertex operator algebras.  We know that this DHR-endomorphism
decomposes into a sum of two irreducible DHR-endomorphisms
$\be_1$ and $\be_2$, thus one of them must have all the vectors
having integer weights and the other has all the vectors having
half-integer weights.   The former  DHR-endomorphism then has
a conformal spin $1$ and the other has $-1$, as desired.

We have two possible fusion rules here, that is, 
$\Z_2\times \Z_2$ and $\Z_4$.
By Rehren \cite{R1}, we have $\om_{\a^n}=(\om_\a)^{n^2}$,
which shows that the fusion rules must be given by 
$\Z_2\times \Z_2$, not $\Z_4$.
\end{proof}

We may and do assume that the DHR sector $\be_1$ has a conformal
spin $1$ and $\be_2$ has $-1$.  By Lemma \ref{spin1}, we
can make an extension of $\A_L^\si$ of index 2 using
the vacuum sector and the sector $\be_1$.  We denote this net
by $\ti\A_L$.  This corresponds to the twisted orbifold
vertex operator algebra $\ti V_L$.  We now
have the following theorem.

\begin{theorem}\label{twisted}
The local conformal net $\ti\A_L$ constructed above is holomorphic
and its vacuum character is the same as the one for the vertex
operator algebra $\ti V_L$.
\end{theorem}

\begin{proof}
By \cite[Proposition 24]{KLM}, we know that the $\mu$-index  of 
$\ti\A_L$ is 1.  By the description of the weights in the
proof of Lemma \ref{L6} above, the decomposition of $\ti\A_L$
for the restriction to $\Vir_{1/2}^{\otimes 2d}$ exactly
corresponds to the decomposition in \cite[Theorem 4.10]{DGH},
in a similar way to the decomposition in Lemma \ref{L6}.
Since the characters for the irreducible DHR endomorphisms
of the Virasoro net $\Vir_{1/2}$ and the characters for
the modules of the Virasoro vertex operator algebra $L(1/2,0)$
coincide,  the vacuum characters for $\ti A_L$ and $\ti V_L$
also coincide.
\end{proof}

In \cite[page 625]{FG}, Fr\"ohlich and Gabbiani have conjectured
that the $S$-matrix arising from a braiding on a system of
irreducible DHR sectors as in \cite{R1}
coincides with the $S$-matrix arising
from the change of variables for the characters.  
(We have called a completely rational net satisfying this
property {\sl modular} in \cite{KL3}.)  If a local
conformal net is holomorphic, this conjecture would imply that
the vacuum character is invariant under the change of 
variables given by $SL(2,\Z)$, that is, a modular invariant
function.  The next corollary shows that this modular invariance
is true for our current examples.

\begin{corollary}
Let $C$ be self-dual double-even binary code of length $d$
and $L$ be $L_C$ or $\ti L_C$, the lattice arising from  $C$ as above.
Let $\A$ be the local conformal net $\A_L$ or its twisted orbifold
as above.  Then the vacuum character of the net $\A$ is invariant
under $SL(2,\Z)$.
\end{corollary}

\begin{proof}
Such modular invariance for the vacuum characters of the vertex
operator algebras $V_L$ and $\ti V_L$ have been proved by
Zhu \cite{Z}.  Thus Theorem \ref{twisted} gives the conclusion.

The conclusion also follows from Proposition 2 in \cite{KL3}.
\end{proof}

\begin{example}\label{moon}{\rm
As in Example \ref{golay}, we consider
the Leech lattice $\La$ as $L$ above.  Then the corresponding
net $\ti A_\La$ has the vacuum character as follows, which is the
modular invariant $J$-function.
$$q^{-1}+196884q+21493760q^2+8642909970 q^3+\cdots.$$
This net $\ti\A_\La$ is the counterpart of the moonshine vertex
operator algebra, and has the same vacuum character as the moonshine
vertex operator algebra.  We call this the {\sl moonshine net} and
denote it by $\A^\natural$.
The near-coincidence of the coefficient
$195884$ in this modular function and the dimension $196883$ of the
smallest non-trivial irreducible representation of the Monster
group, noticed by J. McKay, was the starting point of the
Moonshine conjecture \cite{CN}.  Note that the nets $\A_\La$
and $\ti\A_\La$ both have central charge $24$ and are holomorphic,
but these nets are not isomorphic because the vacuum characters
are different.  To the best knowledge of the authors, this pair
is the first such example.
}\end{example}

\section{Preliminaries on automorphism of framed nets}

For a local conformal net $\A$ on a Hilbert space $H$
with a vacuum vector $\Omega$, its
automorphism group $\Aut(\A)$ is defined to be set of unitary operators
$U$ on $H$ satisfying $U\Omega=\Omega$ and
$U \A(I) U^*=\A(I)$ for all intervals $I$
on the circle.  This group is also called the {\sl gauge group} 
of a net.  (See \cite[Definition 3.1]{X3}, where Xu studies
an orbifold construction with a finite subgroup of $\Aut(\A)$.)

\begin{definition}\label{framed-net}{\rm 
If a local conformal net $\A$ is an irreducible extension of
$\Vir_{1/2}^{\otimes d}$ for some $d$, we say that $\A$ is
a {\sl framed net}.}
\end{definition}

We prepare general statements on the automorphism
groups of framed nets here.
We first recall some facts on automorphisms of a vertex operator
algebra.  Our net $\ti \A_L$ is defined as a simple
current extension by the group action of $\Z_2$, so it has
an obvious dual action of $\Z_2$, and it gives an order 2
automorphism of the net $\ti \A_L$, whose counterpart for
the moonshine vertex operator algebra has been studied by
Huang \cite{Hu}.  It is known that the Monster group has 
exactly two conjugacy classes of order 2, and they are called
2A and 2B.  The above order 2 automorphism of the moonshine
vertex operator algebra belongs to the 2B-conjugacy class.
In theory of vertex operator algebras, a certain construction
of an automorphism of order 2 of a vertex operator algebra has
been well studied by Miyamoto \cite{Mi1}.  In the case of
the moonshine vertex operator algebra, it is known that this
gives an automorphism belonging to the 2A-conjugacy class.
We construct its counterpart for our framed net as follows.

We start with a general situation where we have three
completely rational nets $\A, \B,\C$ with an irreducible
inclusion $\B\otimes\C\subset\A$ and $\B$ is isomorphic to
$\Vir_{1/2}$.  (Later we will take $\A$ to be a framed net.)
Then the inclusion $\B\otimes\C\subset\A$ automatically has
a finite index as in \cite[Proposition 2.3]{KL1}.  Let $\th$
be the dual canonical endomorphism for the inclusion
$\B\otimes\C\subset\A$.  Then we have a decomposition
$\th=\bigoplus_{h,j} n_{h,j} \la_h\otimes\ti \la_j$, where
$\{\la_0=\id,\la_{1/16},\la_{1/2}\}$ is the system of
irreducible DHR sectors of $\B\isom\Vir_{1/2}$,
$\{\ti\la_0=\id,\ti\la_1,\cdots,\ti\la_m\}$ is the system of
irreducible DHR sectors of $\C$, and $n_{h,j}$ is a
multiplicity.  Set $\hat\th=\bigoplus_{h=0,1/2}
n_{h,j} \la_h\otimes\ti \la_j$.  We assume that $\hat\th$ is
different from $\th$, that is, we have some $j$ for which
the multiplicity $n_{1/16,j}$ is nonzero.  Then
by \cite[Corollary 3.10]{ILP}, this $\hat\th$ is a dual
canonical endomorphism of some intermediate net $\hat \A$
between $\B\otimes\C$ and $\A$.  (Note that $\hat\A$ is
local because $\A$ is local.)  We have the following lemma.

\begin{lemma}\label{ind2}
The index $[\A:\hat\A]$ is $2$.
\end{lemma}

\begin{proof}
Localize $\th$ on an interval and let
$N=\B(I)\otimes\C(I)$, $\hat M=\hat\A(I)$, $M=\A(I)$.
The following map $\pi$ gives an automorphism of order
2 of the fusion rule algebra.
\begin{eqnarray*}
\pi(\la_h\otimes\ti\la_j)&=&\hphantom{-}
\la_h\otimes\ti\la_j,\quad{\rm if\ } h=0,1/2,\\
\pi(\la_h\otimes\ti\la_j)&=&
-\la_h\otimes\ti\la_j,\quad{\rm if\ } h=1/16.
\end{eqnarray*}
By the crossed product type description of the algebra $M$ in
\cite[Section 3]{ILP}, this $\pi$ induces an automorphism
of order 2 ot the factor $M$.  (By the assumption on $\hat\th$,
this map $\pi$ is not identity.)
We denote this automorphism again by $\pi$.  It is
easy to see that $\hat M$ is the fixed point algebra of $\pi$ on
$\M$.  Since this fixed point algebra is a factor, $\pi$ is
outer and thus the index $[\A:\hat\A]$ is $2$.
\end{proof}

Thus the extension $\hat\A\subset \A$ is a simple current extension
of index 2, that is, given as the crossed product by an action 
of $\Z_2$, and thus, we have the dual action on $\A$.  Note that
the vacuum Hilbert space of $\A$ decomposes into a direct sum
according to the decomposition of $\th$, and on each Hilbert subspace,
this automorphism acts as a multiplication by 1 [resp. $-1$],
if the Hilbert subspace corresponds to a sector containing
$\la_0,\la_{1/2}$ [resp. $\la_{1/16}$] as a tensoring factor.  So
this automorphism of order 2 corresponds to the involution of a
vertex operator algebra studied by Miyamoto \cite[Theorem 4.6]{Mi1}.
We denote this automorphism by $\tau_B$.  When $\A$ is a framed net,
an irreducible extension of $\Vir_{1/2}^{\otimes d}$, then we can
take the $k$-th component of $\Vir_{1/2}^{\otimes d}$ as above $\B$,
and we obtain an automorphism $\tau_k$ for this choice of $\B$.

Suppose now that the above $\th'$ is equal to $\th$.  That is,
no subsector of $\th$ contains $\la_{1/16}$ as a tensoring factor.
Then we can define another map $\pi'$ as follows.

\begin{eqnarray*}
\pi'(\la_h\otimes\ti\la_j)&=&\hphantom{-}
\la_h\otimes\ti\la_j,\quad{\rm if\ } h=0,\\
\pi'(\la_h\otimes\ti\la_j)&=&
-\la_h\otimes\ti\la_j,\quad{\rm if\ } h=1/2.
\end{eqnarray*}

In a very similar way to the above, this map also induces an
automorphism of order 2 of the net $\A$.  This is a counterpart
of another involution given in \cite[Theorem 4.8]{Mi1}.

We now have the following general description of a framed net.

\begin{theorem}
Let $\A$ be a framed net, an irreducible extension
of $\Vir_{1/2}^{\otimes d}$.  Then there exist integers $k,l$ 
and actions of $\Z_2^k$, $\Z_2^l$
such that $\A$ is isomorphic to a simple current extension
of a simple current extension of $\Vir_{1/2}^{\otimes d}$
as follows.
$$\A\cong (\Vir_{1/2}^{\otimes d} \rtimes \Z_2^k)
\rtimes \Z_2^l.$$
\end{theorem}

\begin{proof}
For the fixed subnet $\B$, consider the group of automorphisms
of $\A$ generated by $\tau_j$, $j=1,2,\dots,d$.  This is
an abelian group and any element in this group has order at most
2, thus this group is isomorphic to $\Z_2^l$ for some $l$.
Note that the net $\A$ is realized as an extension of this fixed
point net by the group $\Z_2^l$ and the dual action of this
group action recovers the original action of $\Z_2^l$ on $\A$.

Consider the fixed point net of $\A$ under this group $\Z_2^l$.
This net contains $\B$ as a subnet, so decompose the vacuum 
representation of this net as a representation of $\B$, then
we have only tensor products of $\lambda_0$, $\lambda_{1/2}$
and each tensor product has multiplicity 1 since each has
dimension 1.  This is then a crossed product extension of
the net $\B$ and the group involved is an abelian group whose
elements have order at most 2, so the group is isomorphic to
$\Z_2^k$.
\end{proof}

\begin{remark}\label{unique-Q}{\rm
When we make a simple current extension of a local conformal net
$\A$ by a group $G$, the extended local net is unique, as long
as the group $G$ is fixed as a subcategory of the
representation category of $\A$.  That is, we have a choice of
a representative, up to unitary, of an automorphism for each
element of $G$, but such choices produce the unique extension.
(In general, we may have two $Q$-systems depending on the choices,
and they may differ by a 2-cocycle $(c_{gh})_{gh}$
on $G$ as in \cite{IK}, but locality forces $c_{gh}=c_{hg}$, which
in turn gives triviality of the cocycle.  So all the local $Q$-systems
realizing $G$ are unitarily equivalent.)  We have a two-step extension
of $\Vir_{1/2}^{\otimes d}$, thus the $Q$-system of this extension
and the resulting extended net are both
unique as long as the two groups $\Z_2^l$, $\Z_2^k$ are fixed in
the representation categories.
}\end{remark}

\section{The automorphism group of the moonshine net}

In this section, we make a more detailed study of the automorphism
group and show that it is indeed the monster group. 

We begin with some comments.

Starting with a M\"obius covariant local family of Wightman fields
$\{\phi_a(f)\}$, one can (obviously) canonically construct a M\"obius
covariant local net $\A$ if each operator $\phi_a(f)$ is essentially
selfadjoint on the common invariant domain and bounded continuous
functions of $\phi_a(f_1)$ and $\phi_b(f_2)$ commute when the real
test functions $f_1, f_2$ have support in disjoint intervals of $S^1$. 
In this case we shall say that the family $\{\phi_a(f)\}$ is
\emph{strongly local}.

With $V$ a vertex algebras equipped with a natural scalar product, we
shall denote by $\Aut(V)$ the automorphism group of $V$ that we define
as follows: $g\in\Aut(V)$ if $g$ is a linear invertible map from $V$
onto $V$ preserving the inner product with $gL = Lg$ for any M\"obius
group generator $L$, $g\Omega = \Omega$, and such that $gY(a,z)g^{-1}$
is a vertex operator of $V$ if $Y(a,z)$ is a vertex operator of $V$. 
(One does not have a positive definite inner product for a general
vertex operator algebra, but we are to consider the moonshine vertex
operator algebra, and this has a natural scalar product preserved by
the action of the Monster group by \cite[Section 12.5]{FLM}, so we can
use this definition of the automorphism group for a vertex operator
algebra.)  By the state-field correspondence it then follows that
$gY(a,z)g^{-1}=Y(ga,z)$.  To simplify the notation we shall denote by
$g$ also the closure of $g$ on the Hilbert space completion of $V$,
which is a unitary operator.

Let now $L(1/2,0)^{\otimes 48}\subset V^\natural$ be the realization
of the moonshine vertex operator algebra as a framed vertex operator
algebra.  Let $\Vir_{1/2}^{\otimes 48} \subset \A^\natural$ be the
corresponding realization of the moonshine net as a framed net as in
Example \ref{moon}.  (By the ``correspondence'', we mean that the two
$Q$-systems are identified.)  Note in particular that $V^\natural$ is
equipped with a natural scalar product and the completion of
$V^\natural$ under with respect to this inner product is denoted by
$H$.

\begin{lemma}\label{gen}
The vertex operator sub-algebras $V$ of $V^\natural$ generated by
$g(L(1/2,0)^{\otimes 48})$ for all $g\in\Aut(V^\natural)$ is equal to
$V^\natural$.
\end{lemma}

\proof
The subspace $V_2$ contains the Virasoro element of $V^\natural$ and some
other vector which is not a scalar multiple of the Virasoro
element, and it is a subspace of $V^\natural_2$, 
and has a representation of $\Aut(V^\natural)$, the monster group.
The representation of the monster group on the space $V^\natural_2$
is a direct sum of the trivial representation and the one
having the smallest dimension, 196883, among the nontrivial
irreducible representations.  From these, we conclude
$V_2=V^\natural_2$.  Then by \cite[Theorem 12.3.1 (g)]{FLM}, we
know that $V=V^\natural$.
\endproof

Let $H$ be the Hilbert space arising as the completion of $V^\natural$
with respect to the natural inner product.  Because of the
identification of the $Q$-systems, we may regard that the net
$\A^\natural$ acts on this $H$.

Denote $\B$ the local conformal net associated with $L(1/2,0)^{\otimes 48}$.
Clearly $\B = \bigotimes_{k=1}^{48}\B_k$ where the $\B_k$'s are 
all equal to the local conformal net $\Vir_{1/2}$.

Let $T_k(z)$ and $T(z)$ denote the stress-energy tensor of $\B_k$ and
$\B$.  We consider $T_k(z)$ as acting on the Hilbert space $H$ of $\B$
(correspondingly to the obvious identification of of the Hilbert space
of $\B_k$ as a subspace of $H$).  As is well known each $T_k$ is a
strongly local Wightman field on the domain of finite-energy vectors
for $T_k$ (which is dense on the Hilbert space of $\B_k$) \cite{BS}. 
Therefore, by tensor independence, the full family
$\{T_k\}_{k=1}^{48}$ (and $T$) is a strongly local family of Wightman
fields on $H$.  Furthermore we have:

\begin{lemma} The family of Wightman fields $\{T_{k,g}\equiv gT_k
g^{-1}: g\in \Aut(V^\natural), k= 1, 2,\dots 48\}$ is strongly local.
\end{lemma}

\proof
Fix $g,g'\in\Aut(V^\natural)$ and $k,k'\in\{1,2,\dots 48\}$. As vertex 
operators (formal distributions) we have
\[
T_{g,k}(z)T_{g',k'}(w) = T_{g',k'}(w)T_{g,k}(z), \quad z\neq w
\]
because $T_{g,k}$ and $T_{g',k'}$ belong to $V^\natural$.
Therefore
\[
T_{g,k}(f)T_{g',k'}(f')v = T_{g',k'}(f')T_{g,k}(f)v, 
\]
if $f,f'$ are $C^{\infty}$ functions with support in disjoint 
intervals of $S^1$ and $v$ is a vector of $V^\natural$.

Recall that the vectors of $V^\natural$ are the finite energy vectors 
of $L_0$, the conformal Hamiltonian of $T$. Denoting by $L_{0;k,g}$ 
the conformal Hamitonian of $T_{g,k}$ we have the energy bounds 
in \cite{BS}
\[\label{bound}
||T_{g,k}(f)v||\leq c_{f}||(L_{0;k,g} + 1)v ||\leq c_{f}||(L_0 + 
1)v||
\]
where $v$ is a finite energy vector, $c_f$ is a constant, and in the
last inequality we have used that $gTg^{-1} = T$ and
\[
L_0 = \sum_k L_{0;k,\iota}= \sum_k L_{0;k,g}
\]
(on finite-energy vectors) thus $L_{0;k,g}< L_0$.

Analogously, as in \cite{BS}, one can verify the other bound necessary 
to apply a theorem in \cite{DF}
\[
||[T_{g,k}(f),L_0]v||\leq c_{f}||(L_0 + 1)v||\ ,
\]
check the required core conditions, and conclude that the $T_{g,k}$'s form 
a strongly local family of Wightman fields.
\endproof

\begin{corollary}
Let $\tilde\A$ be the conformal net on $H$ generated by all these
fields $T_{g,k}$, namely by all the nets $g\B_k g^{-1}$,
$g\in\Aut(V^\natural)$ and $k\in\{1,2,\dots 48\}$.  Then $\tilde\A$ is
\emph{local} and we may identify $\tilde\A$ and $\A^\natural$.
\end{corollary}

\proof The locality of $\tilde\A$ is is immediate from the above
lemma.  We will only need to check the cyclicity of the vacuum 
vector $\Omega$.

Let $E$ be the selfadjoint projection of $H$ onto the closure of
$\ti\A\Omega$.  Then $E$ commutes with $\ti\A$ and with the unitary
M\"obius group action.  The unitary $U\equiv 2E - 1$ thus fixes the
vacuum, maps finite energy vectors to finite energy vectors and $U
T_{g,k} U^* = T_{g,k}$.  By Lemma \ref{gen} $UY(a,z)U^* = Y(a,z)$ for
every vertex operator $Y(a,z)$ of $V^{\natural}$.  Thus $Ua =
UY(0,a)\Omega =Y(0,a)U\Omega = Y(0,a)\Omega = a$ for every vector 
$a\in V^{\natural}$, namely $U$ and $E$ are the identity, so $\Omega$ 
is cyclic.

Last, the identification of $\tilde\A$ and $\A^\natural$ follows by
Remark \ref{unique-Q}.  
\endproof

\begin{theorem}\label{Aut}
There is a natural identification between $\Aut(V^\natural)$, the
automorphism group of the moonshine vertex operator algebra, and
$\Aut(\A^\natural)$, the automorphism group of the moonshine local
conformal net.  Thus $\Aut(\A^\natural)$ is the Monster group by
\cite{FLM}.
\end{theorem}

\proof 
Since $\A^\natural(I)$ is generated by $g\B(I)g^{-1}$ as $g$ varies in
$\Aut(V^\natural)$, clearly every $g\in\Aut(V^\natural)$ gives rise to 
an automorphism of $\A^\natural$.

It remains to prove the converse.  So let $g\in\Aut(\A^\natural)$ be
given; we shall prove that $g$ corresponds to an automorphism of
$V^\natural$.

By definition $g$ is a unitary operator on $H$ that implements an
automorphism of each $\A^\natural(I)$ and fixes the vacuum.  Then $g$
commutes with the M\"obius group unitary action \cite{BGL}, in particular with
the conformal Hamiltonian, (the M\"obius group action is determined by
the modular structure); therefore $g$ maps finite energy vectors to
finite energy vectors, i.e. $gV^\natural = V^\natural$.

Since $V^\natural$ is generated by the vertex operators $T_{k,g'}$ 
(as $g'$ varies in $\Aut(V^\natural)$), 
and $\Aut(\A^\natural)\supset \Aut(V^\natural)$, it is sufficient 
to prove that $gT_k g^{-1}$ is a vertex operator of $V^\natural$.

The proof is based on the state-field correspondence. Set $W(z) = 
gT_k(z)g^{-1}$ and recall that $W$ is a Wightman field, in particular 
$W(z)\Omega$ is an analytic vector valued function for $|z|<1$. 

With $f$ a real test function with support in an interval $I$, the
operator $W(f)$ is affiliated to the von Neumann algebra
$g\A^\natural(I)g^{-1} = \A^\natural(I)$ and it thus follows that $W$
is strongly local with respect to $\{T_{k,g'}\}$,
$g'\in\Aut(V^\natural)$ and $k\in\{1,2,\dots 48\}$.  Since these
fields generate $V^{\natural}$, by Dong's lemma \cite{Kac} $W$ is
local with respect to all vertex operators of $V^{\natural}$.

Set $a\equiv W(0)\Omega$ and let $Y(a,z)$ be the vertex
operator of $V^\natural$ corresponding to $a$ in the state field
correspondence, thus $Y(a,0)\Omega = a$. 

We have
\[
W(z)\Omega = e^{zL_{-1}}W(0)\Omega = e^{zL_{-1}}a = 
e^{zL_{-1}}Y(a,0)\Omega =Y(a,z)\Omega, \quad |z|< 1,
\] 
where $L_{-1}$ is the infinitesimal translation operator of
$V^{\natural}$.  Now we have for $z\neq w$ and $b\in V^{\natural}$
(cf.  \cite[Sect.  4.4]{Kac})
\[
W(z)e^{wL_{-1}}b\Omega = W(z)Y(b,w)\Omega = 
Y(b,w)W(z)\Omega = Y(b,w)e^{zL_{-1}}a\Omega ;
\]
here we have made use of the mutual locality of $W(z)$ and $Y(a,w)$
and that the infinitesimal translation operator of $W$ is $L_{-1}$
because $gL_{-1}g^{-1} = L_{-1}$. Letting $w=0$ a we get
\[
W(z)b = Y(b,0)e^{zL_{-1}}a\ .
\]
Analogously $Y(a,z)b = Y(b,0)e^{zL_{-1}}a$, so $W(z)= Y(a,z)$ as desired.

\medskip

\noindent{\bf Acknowledgments.} A part of this work was done during
visits of the first-named author to Universit\`a di Roma ``Tor
Vergata'' and a part of this paper was written during a visit of the
first-named author at the Erwin Schr\"odinger Institute.  We
gratefully acknowledge the support of GNAMPA-INDAM and MIUR (Italy)
and Grants-in-Aid for Scientific Research, JSPS (Japan).  The
first-named author thanks T. Abe, A. Matsuo, M. Miyamoto, H. Shimakura
and H. Yamauchi for helpful discussions on vertex operator algebras.
We thank M. M\"uger for showing us the unpublished thesis \cite{S} of
Staszkiewicz and explaining its contents, of which we were first unaware.
We also thank B. Schroer and F. Xu for their comments on the first
version of this paper.
Finally we are grateful to S. Carpi and M. Weiner for helpful
discussions concerning Theorem \ref{Aut} (also as part of a joint
work in progress).

\end{document}